\documentclass{amsart}

  \usepackage[all]{xy}

  \usepackage{epsf,epsfig,amsfonts,graphicx,color}

\usepackage{hyperref}
\hypersetup{draft=false}

\hypersetup{citecolor=DeepPink4}
\hypersetup{linkcolor=DarkRed}
\hypersetup{urlcolor=DarkBlue}

\usepackage{cleveref}

\usepackage{amsmath,amssymb}
  \usepackage{url}
   
 \def\R{\mathbb R}

\def\C{\mathbb C}

\newcommand{\beq}{\begin{equation}}
\newcommand{\eeq}{\end{equation}}

\begin{document}


\author{Richard Montgomery}
\address{Mathematics Department\\ University of California, Santa Cruz\\
Santa Cruz CA 95064}
\email{rmont@ucsc.edu} 

\title{Dropping  Bodies} 
\thanks{ }

\maketitle

  Drop three  bodies.  Where can they go?  
   
  When  someone says they've dropped a cup, we  imagine it falling down to the ground, attracted to the earth by gravity.   But please, in this  thought experiment,  take away the earth.  There's no direction ``down''.  Our three bodies are alone in the universe,   attracted only to each other.  
   
Take each body  to be  a point mass.   Dropping  the bodies means letting  them    go from rest,    subject only to the rules of  Newtonian mechanics 
and     the assumption that the  only   forces acting on a body are   the   gravitational  $1/r^2$ pulls of the other two.   
 Each body  will  then sweep out a plane curve, subject to the attractive pull of the other two moving bodies.
   Taken together, these three parameterized plane curves form
 a  solution to   the Newtonian three-body problem.
 
\begin{figure}
  \includegraphics[width=6cm]{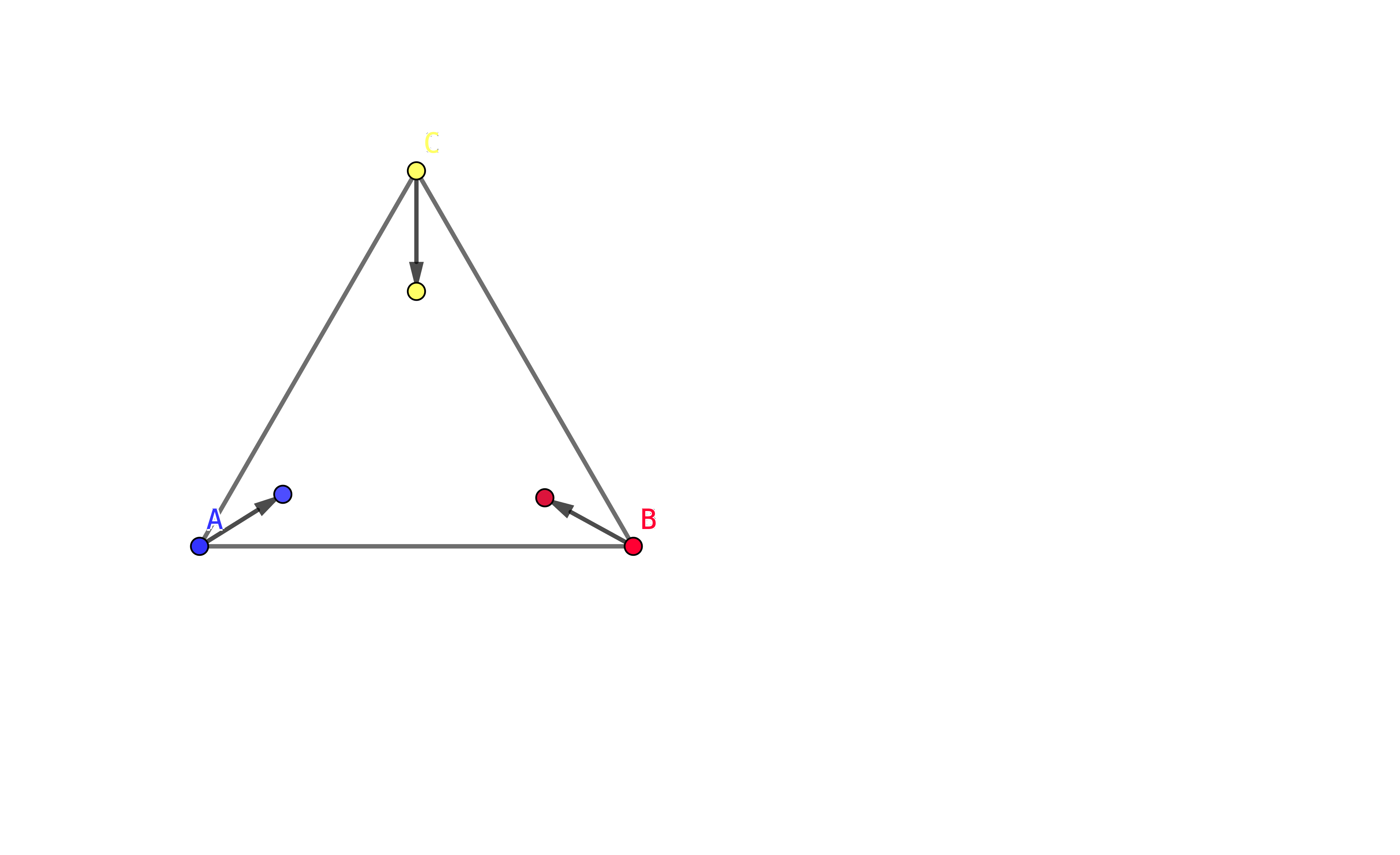}
 \caption{Dropping an equilateral triangle.  Lagrange showed the triangle shrinks
 to triple collision, remaining equilateral at each instant.} 
 \label{figLagrange}
 \end{figure}

\begin{figure}
  \includegraphics[width=6cm]{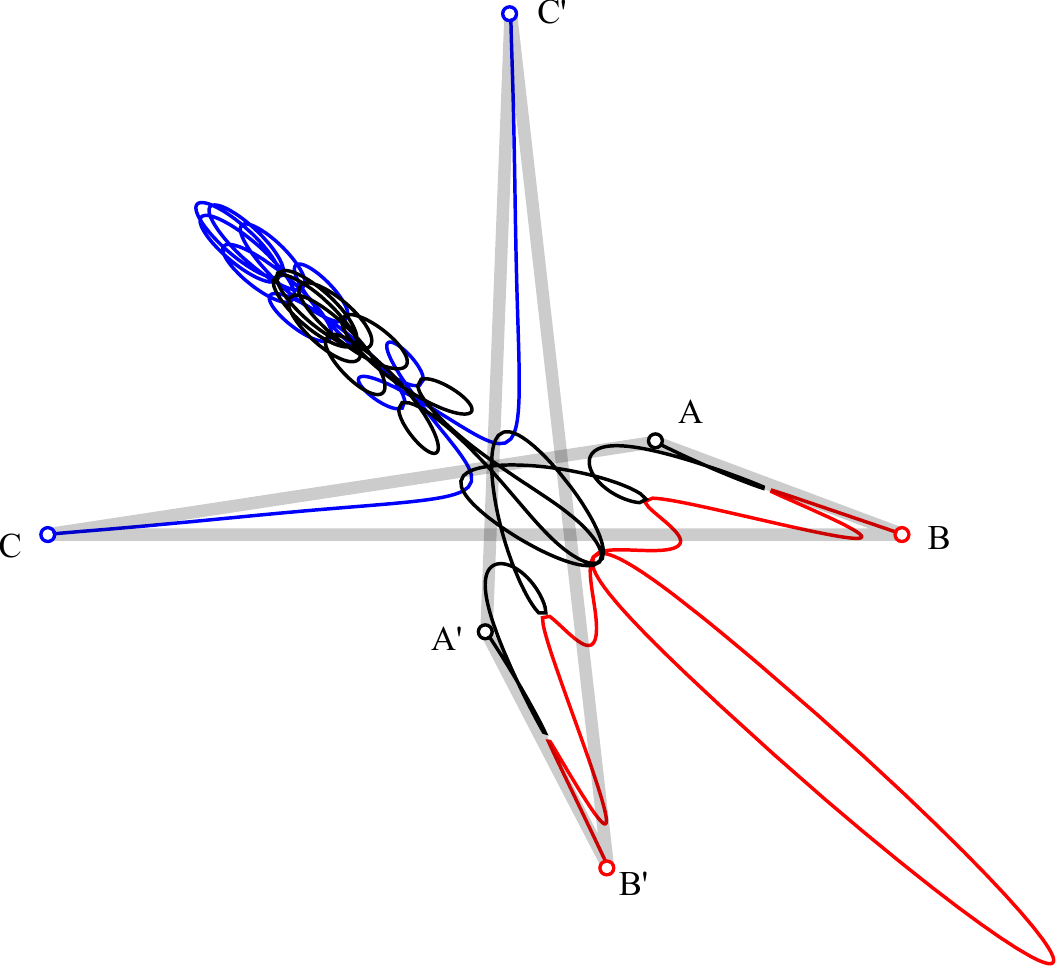}
 \caption{Orbit 18  or ``$F_{18}(1,1,1)$'' in   the Li-Liao database \cite{Li}.    
 Drop triangle A B C.  A half-period later we will arrive  at its reflection,  triangle A' B' C' .
 The two labelled  triangles are related by a reflection.  As a consequence,  at the 
 half-way time between, the three bodies must  line  up along  the reflectional line of symmetry of the figure.  See the
 section `Symmetry Puzzles''. } 
 \label{fig18Gofen}
 \end{figure}

\begin{figure}
  \includegraphics[width=6cm]{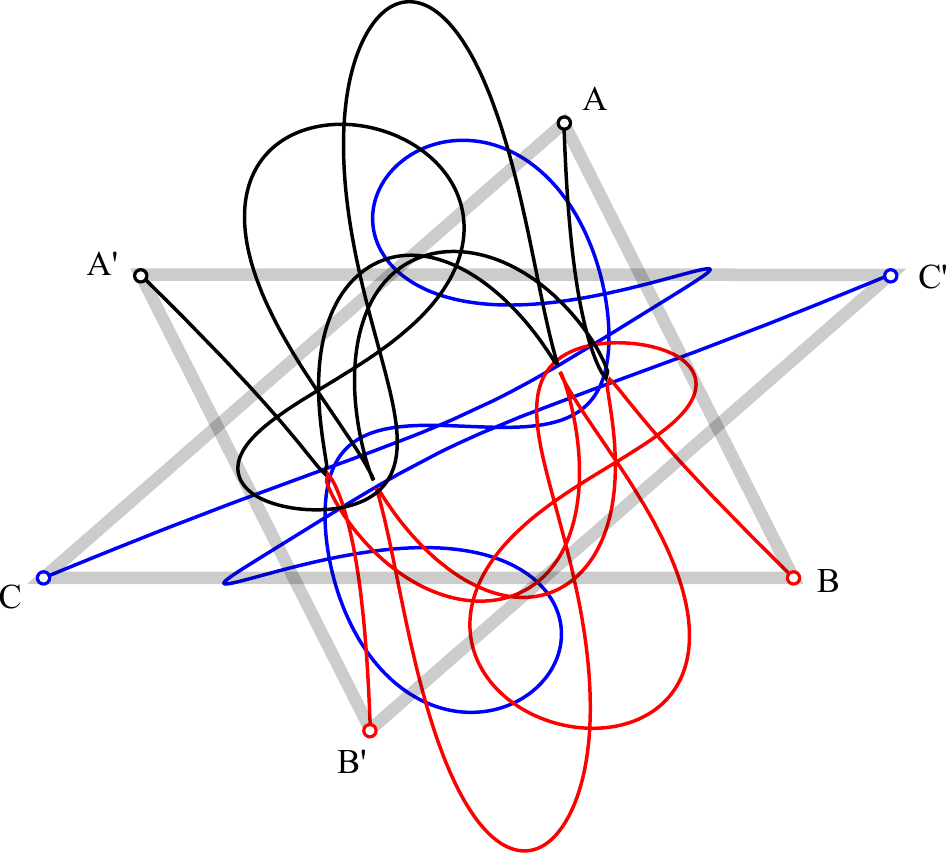}
 \caption{Orbit 14   or ``$F_{14}(1,1,1)$'' in  the Li-Liao database \cite{Li}.  The two {\it unlabelled} triangles are related by  central inversion, i.e. a 180 degree rotation.
 That inversion takes $C$ to $B'$ rather than $C'$.}
 \label{fig14Gofen} 
 \end{figure}

\begin{figure}
  \includegraphics[width=6cm]{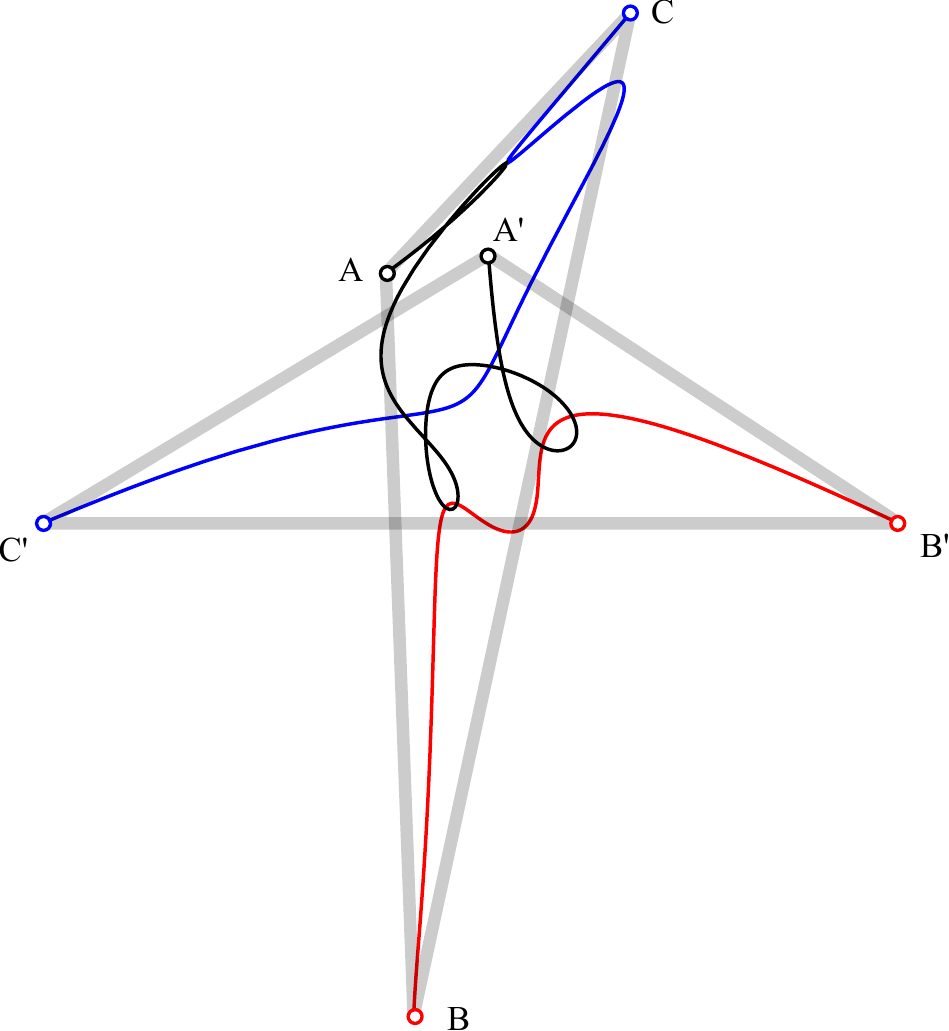}
 \caption{Orbit 1 or ``$F_{1}(1,1,1)$'' in   the Li-Liao database \cite{Li}. 
  No isometries relate the two braked triangles.}
 \label{fig6Gofen}  
 \end{figure}

\begin{figure}
  \includegraphics[width=6cm]{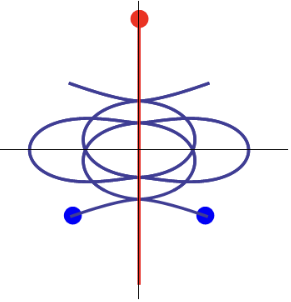}
 \caption{A periodic brake orbit in which the three bodies form an isosceles triangle
 at each instant.  The orbit  discovered by Nai Chia Chen attempts to redraw an early Paul Klee painting.} 
 \label{fig1NaiChia}  
 \end{figure}

\begin{figure}
  \includegraphics[width=6cm]{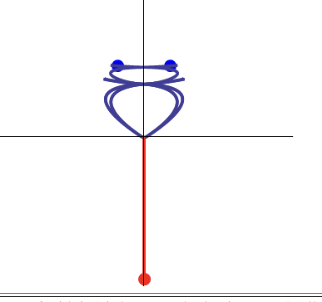}
 \caption{Another one of the $6 * \infty$ periodic brake orbits found by Chen.  } 
 \label{fig2NaiChia}  
 \end{figure}

Below we have depicted  a half-dozen  answers  to our question - indications of where the bodies went.   
    The answers vary widely depending on the starting triangle.   
  We urge the  reader to view some of the  animations  which can be found at \cite{Gofen} and \cite{Chen}.
     In the sampled pictures and animations the  three masses   are   equal.  
 (The answers   seem to be prettier that way.)  The solution depicted in figure  \ref{figLagrange} was found by Lagrange.
         The other   figures depict periodic   solutions which are  solutions for which the three bodies shuttle  back and forth between
       two  ``brake triangles''-  configurations of the three bodies at which they are  instantaneously at rest or ``braked''. Either 
       brake triangle  can be supposed to be  the initial configuration from which the bodies are dropped. 
Figures  \ref{fig18Gofen}, \ref{fig14Gofen}, and \ref{fig6Gofen}   have been selected  from a database of thirty  collision-free equal mass `dropped' solutions  found by  Li and Liao (\cite{Liao}) and available to view at (\cite{Li}).  They are all
collision-free.  Figures \ref{fig1NaiChia}, and \ref{fig2NaiChia} display two of the infinitely many periodic
  isosceles  brake  solutions   found by  Nai-Chia Chen (\cite{Chen2}), all of which suffer binary collisions
and whose existence she established  in her PhD thesis.  
Chen's discoveries remind me of Paul Klee's drawings.  Her  animations are as if Klee's
 drawings had children with  the mobiles of   Alexander Calder.  
 
Just above, we used the term   brake orbit as a synonym for dropped solution.   A brake orbit is thus a  solution to the N-body problem
 for which, at some instant, all N bodies are instantaneously at rest, or in
 automotive parlance, have ``braked'' to a stop.   Another synonym for brake orbit is ``free-fall solution''. 
  
   The  observer  will
    see various  symmetries within the orbits.   For example,   figure (\ref{fig18Gofen}) 
  has a reflectional symmetry.  If you watch the animation  
  that reflectional  symmetry  become spatio-temporal.   When we drop the bodies at time $t = 0$ and if
  $T$ is the period, then the other brake triangle arises at time $T/2$.  At time   $T/4$   all three bodies    lie on the  symmetry line of the figure.
  Why? 
  

\section{History}  Lagrange showed  that a dropped equilateral triangle   stays equilateral all the way  until it
has collapsed to triple collision.   This is not  a surprise when the masses are equal
but his result holds true for unequal masses.    Euler, five years before Lagrange, 
 uncovered a similar fact by placing    three masses judiciously on a line:
he could get their shape to  remain  the same upon dropping them.   It is easy to get the placement right for three equal masses. Place one at the midpoint of the
other two. The center of mass of the triple is then   this midpoint and the extremal  two  collapse symmetrically
 onto this midpoint.
For general unequal masses Euler needed  to solve a quintic  to get the mass placements right.

In the 1893 a person named Meissel asserted, with little evidence,  that if   masses   in the ratio of $3:4:5$ are
  placed at the vertices of a Pythagorean $3:4:5$ triangle, with mass $5$ placed at the vertex opposite the side of length $5$
  and so on, and  the bodies are dropped, then the resulting solution is periodic.    Verifying Meissel's assertion became a  benchmark 
   problem for numerical computation in celestial mechanics and was christened   the ``Pythagorean three-body problem'' or  ``Burrau's three-body problem''.
 Carl Burrau's name (\cite{Burrau}) became attached to the problem after  he  took up Meissel's challenge and published the inconclusive results of 
hand-cranked  numerical integration for the problem in 1913.  Let's all  raise our glasses to the  long-suffering degree candidate
Sigurd Kristensen  who, 109 years ago,  did   an essential part -``Ein wesentlicher Teil der Rechernarbeit '' of Burrau's computational work  (one guesses all)  
but whose name has faded from view,     no co-authorship,    only an honorable mention within the article.   In 1967 the celestial
  mechanician Szebehely   and his team  proved Meissel  wrong using the new-fangled digital computers. They showed
  that  the smaller of the two  masses eventually form  a binary pair which escapes to   infinity.  Before the big  escape many interchanges and
  several  near collision incidents
  occur, incidents requiring  Szebehely and team  to  implement   Levi-Civita's regularization of binary collisions in order to get an accurate
  integration and proceed to follow the orbit to its end \footnote{Nowadays, high order Taylor methods can 
  recreate Szebehely's discovery without requiring  
 regularization.  See \cite{Gofen} and the final section of this paper.} 
 
   Wait! What  is Levi-Civita doing  here?  Didn't he work on connections and Riemannian geometry and write a book with his advisor Ricci
   on tensor calculus?   A few years after Levi-Civita 
  uncovered the  famous connection now bearing his name, 
  he established a suprisingly  simple change of variables which ``regularizes'' the binary collisions
  of the N-body problem.  Gravitational forces blow up at collisions, making it seem  impossible to continue
  solutions through collision.   However  Levi-Civita's change of both time and space variables   renders the defining ODEs analytic near
  isolated  binary collisions.
  Chen's solutions have binary  collisions and to make sense of them, or  the isosceles three-body problem in general,
  requires  a variant of Levi-Civita's regularization. 
    
  Soon  after Szehebely had established that Meissel was wrong he decided to tweak Meissel's initial conditions
  and found,   remarkably, 
    near to Meissel's initial conditions, a  periodic brake orbit.  At the half-period this solution suffers a binary collision
    in which the non-colliding mass has braked to a stop, and as a consequence of a closer look at
    Levi-Civita, the solution reverses its path, returing to the starting point.   This solution provided inspiration for many.
    At my own university, the astronomer Greg Laughlin (sadly for me, now moved to Yale) collaborated with a dance professor to choreograph 
     Szehebely's near Pythagorean periodic brake orbit.  Three humans played the falling bodies.   See \cite{Laughlin} for more. 
    
  The Newtonian N-body problem is a special case of a  Hamiltonian dynamical system whose
  energy, the Hamiltonian, has the form of kinetic plus potential.   See equation (\ref{N}) below, allowing $V$ to be general.   Brake orbits makes sense in this
  general context.    Seifert, the topologist, wrote a beautiful  paper (\cite{Seifert}) on brake  orbits  in this more general context.
  His paper provided some of the fuel for Weinstein to make a now-rather-famous conjecture known as the Weinstein conjecture - a kind of  contact version
  of the Arnol'd conjectures which drove the development of the field of  symplectic topology.
  The Weinstein conjecture was proved for the case of 3-dimensional contact geometries by Taubes \cite{Taubes} about a decade ago. 
  
  Weinstein's student Ruiz  (\cite{Ruiz1}, \cite{Ruiz2}) coined the name ``brake orbit'' in his thesis   extending Seifert's results.  
  Ruiz insisted   brake orbits be periodic. We prefer to    use the word in our less restrictive sense, requiring only one brake instant while Ruiz
  definition requires two.  Periodic brake orbits must of necessity shuttle back
  and forth between two brake configurations.    
 
 Nothing's sacred here about ``three''.  Drop N bodies. The analogues of
 the solutions found by Euler and Lagrange are nowadays called central configurations.
 In other words, if upon being dropped, the  N bodies maintain their shape while collapsing to total
 collision then that initial configuration is called a ``central configuration''.    If we apply an isometry or a scaling to 
 a central configuration we get another one.   Modulo such isometries and scaling, is it true that the number
 of central configurations is finite?   This problem  made it onto Smale's list of mathematical
 problems for the 21st century (\cite{Smale}),  published in this same esteemed journal.  
 In 2006 (\cite{Hampton})  the problem was finally answered `yes' for the case $N=4$.  
 We have an ``almost yes'' for $N =5$ (see \cite{Albouy}) and a ``we're basically clueless beyond numerical experiments which indicate yes'' for $N > 5$.

\section{The problem}

  The ODE   defining the three-body  problem, some of whose   solutions are  depicted in our figures,   can be written 
   \beq 
   \ddot q = - \nabla V (q).
   \label{N}
   \eeq
   Here $q = (q_1, q_2, q_3)$ records the positions of the three bodies, so that $q_a \in \R^2 \cong \C$,  $a =1,2,3$. 
  The curves $q_a (t)$ are parameterized by Newtonian time $t$.  The double dots over $q$ denotes acceleration -  the second derivative with respect to time.
   $\nabla V$
   denotes the gradient of the potential function (\ref{eq: potential}) with respect to an inner product on the space $\R^2 \times \R^2 \times \R^2$
   of $q$'s called the
   ``mass inner product'':  $$\langle \dot q , \dot q \rangle = \Sigma m_a | \dot q_a |^2.$$
      A dropped solution, or brake orbit,  is one for which, at the initial time $t =0$ we have $\dot q_a = 0, a =1,2, 3$.
   
     The Newtonian potential is
 \beq 
 V = - G \Sigma \frac{m_a m_b}{r_{ab}}, \qquad \text{ where } r_{ab} = |q_a - q_b|
 \label{eq: potential}
 \eeq
is the distance between the bodies. $G >0$ is the gravitational constant, needed if the units on both sides of (\ref{N}) are to match up.
The positive constants $m_a$ are the masses. We set them   equal to each other for the equal mass three-body problem.   

The dynamics defined by equation (\ref{N}) leaves the center-of-mass subspace $\Sigma m_a q_a = 0$ invariant.
This means is that if we start with initial conditions $q_a (0),  \dot q_a (0), a = 1, \ldots , N$ for which  
$\Sigma m_a q_a (0) = 0$ and $\Sigma m_a \dot q_a (0) = 0$ then, for all time the resulting  solution  $q_a (t)$
lies in the subspace:  
$\Sigma m_a q_a (t) = 0$.  There is a standard undergraduate physics trick which reduces any solution to 
such a center-of-mass zero solution. We invoke this reduction to center-of-mass several times   below, perhaps once or twice without saying so.   

 The mass inner product was built so that the usual kinetic energy is written as $K(\dot q) = \frac{1}{2} \langle \dot q , \dot q \rangle$.
 The total energy $K(\dot q) + V(q)$  is conserved, meaning constant along solutions to equation (\ref{N}).
 In addition to this energy, the  linear momentum,
   and angular momentum are conserved.  The linear momentum we've already seen. It is $\Sigma m_a \dot q_a$.
   The angular momentum is given by $\Sigma m_a q_a \wedge \dot q_a $ where
   the ``$\wedge$'' denotes the two-dimensional version of the cross product: $(x, y) \wedge (\dot x, \dot y) = x \dot y - y \dot x$,
   for $(x, y), (\dot x, \dot y) \in \R^2$.    The linear and angular momenta are   linear in velocities, while the energy is
   of the form (kinetic) + (potential) where the kinetic energy is positive definite in velocities and  the potential is an everywhere negative function.  It follows
   that  all our dropped   solutions have   zero linear momentum,  zero angular momentum, and negative energy.

The physically literate reader may have  protested at my potential and said  ``The potential you wrote down is the potential
coming from the fundamental solution of the Laplacian in $\R^3$, not $\R^2$!   It yields the
gravitational  force for  bodies moving in space, not in the plane.  Your $q_a$ must lie  in 
Euclidean space  $\R^3$. Your insistence   $q_a \in \R^2$ is unsightly and   wrong. ''  I  will counter  by
reminding  my literate reader that since  dropped solutions have zero angular momentum 
 they   necessarily remain in the plane containing their initial triangle with vertices  
$q_a (0)$, $a =1,2, 3$   at time $t = 0$. Identify this plane with $\R^2$ and the dynamics are correct.

  \section{Symmetry Puzzles}

  The seed of this paper were two mysteries and  a paradox which Alex Gofen brought to my attention.
  Gofen had been going through a database   of 30 periodic  collision-free brake orbits for the equal-mass three body problem
  which are  compiled in \cite{Liao}.  Gofen was verifying and exploring these solutions 
  using  his own ``Taylor Center''  \cite{Gofen} integration scheme.   A 
  periodic brake orbit   must  shuttle back and forth between two distinct brake configurations, which is to say
  two distinct   brake triangles for  the
  three body problem. Gofen noticed that in 12 out of the 30 cases the  two brake triangles were related by
  an isometry.  Thus mystery number one: why such a large number?  And is the number ``right''?  For example, 
  if we could make a database of the   `next' 300 or 3,000 equal mass periodic collision-free brake orbits
  would we continue to find roughly $1:3$ of them  had congruent brake triangles?  Now to mystery two. 
  One observes in the database  that a symmetry relating  the two  `end'  brake triangles  induces a symmetry of the entire spatio-temporal structure of 
  the orbit.   Why?   An instance of     mystery number two   was discussed above in  relation to figure \ref{fig18Gofen}.
  Finally, to the paradox. Gofen observed that the  solution depicted
  in  figure \ref{fig14Gofen} shuttles back and forth between    two brake  triangles which are  congruent  by
  a rotation, indeed by rotation by 180 degrees.    My own ``shape space'' perspective (see \cite{me_monthly}) on the three-body problem
  told me that what he observed was impossible.   Hence the paradox.

  To  resolve the paradox,   let's begin by understanding  the shuttling back-and-forth.  
  Newton's equations (\ref{N})  enjoy time reversal symmetry: if $q(t)$ solves, so does $q(-t)$.    Now if  $q(t)$ is a brake orbit
  with brake instant $t=0$ then $q(-t)$ has precisely the same initial conditons - same configuration and same zero velocity --
  at time $t=0$ as $q(t)$.
  It follows by the unique dependence of solutions on initial conditions that we must  have  $q(t) = q(-t)$.    Suppose in addition that the brake orbit is periodic with period $T$ .  Then we have the additional temporal symmetry  $q(t+ T) = q(t)$.    By substituting $t = h-T/2$ into this periodicity  relation we get
  $q(h + T/2) = q(h-T/2)$.  Use the time reversal symmetry to get $q(h+T/2) = q(T/2 -h)$.
  Differentiate with respect to $h$ at $h =0$ to see that  $t = T/2$  is  another brake instant.
  Thus periodic brake orbits must shuttle back and forth between two brake configurations.   
  
  This second brake  configuration  must be different from the first. If not, we can cut our period in half and repeat the argument.
  Eventually by this process we either arrive at a fundamental minimal period with two distinct configurations or the
  periods go to zero which means our original brake orbit was in fact a fixed point - a critical point of the potential.
  For  Newton's potential this last possibility is excluded.  The potential  has no critical  points.  N stars cannot just sit in space, attracting each other but not moving.

  Let me proceed now to  shape space thinking.  The Newtonian N-body problem enjoys spatial   symmetries in addition to its temporal symmetries.  
  These are the isometries of space, or, in the case of the planar N-body problem,
  the   isometries of
  the plane.  What this means is that 
   if $R$ is any isometry of the plane and $q(t) = (q_1 (t), q_2 (t), q_3 (t))$ solves the planar three-body problem, then so does $Rq (t)$
  where by $R(q_1, q_2 , q_3)$ we mean $(Rq_1 (t), R q_2 (t), R q_3 (t))$.     We
  can use these symmetries to    push down  the ODEs defining the three-body problem down to a space I call ``shape space''.  
  The points of shape space are {\it oriented} congruence classes of planar triangles.
  Two triangles represent the same ``shape,'' or  oriented congruence class if there is
  an orientation preserving isometry $R$,   i.e. a rotation composed with a translation, taking one to the other.
  
  Here is the salient point of the paradox: when the angular momentum is zero this pushed-down
  dynamics is also of Newtonian type, so the argument of the preceding paragraph holds. And brake orbits all have zero angular momentum.
    A  periodic brake orbit, viewed
  in shape space, is a periodic brake orbit down in shape space, and so  must shuttle back and forth between two shapes and {\it these shapes must be distinct}.
    But Gofen told me he had found orbits where the two brake triangles were related by rotation hence
    down in shape space it was shuttling back and forth between the same point!  
  
  The resolution of this paradox is that Gofen was viewing his  triangles as unlabelled while my triangles have to   be  labelled triangles
  in order for me to  construct shape space with its  dynamics.  
 When  a   rotation $R$ takes a labelled triangle to another
 it must, by definition,   preserve the (mass or vertex) labellings.  But in the example depicted in figure (\ref{fig14Gofen}) the two brake triangles
   are congruent  as {\it unlabelled} triangles, not as labelled one.  
 A 180 degree rotation does take one triangle to the other but in so doing it    messes up their labellings.
  
  When the masses are all equal, the N-body problem enjoys additional   symmetries beyond the Galilean symmetries of time   and
  space isometries. We may  interchange any two masses:  if $(q_1 (t), q_2 (t), q_3 (t))$ solves the equal mass three-body problem then 
  so does $(q_2 (t), q_1 (t), q_3 (t))$, etcetera. 
  The operation of interchanging masses defines a representation of the permutation group on the configuration space
  so that we could write the above interchange of masses 1 and 2 as $\sigma_{12}$.  The brake triangles
    of figure (\ref{fig14Gofen}) are related by a symmetry of the form $F = R \circ \sigma$ where $\sigma$ is one of the transpositions.  
  Although $R$ acts as the identity   on shape space,  such an  $F$ does not.  
  The shape of $q(0)$ and $F(q(0))$ are different, allowing us to avoid  the paradox.
 
  Having extricated ourselves from our mathematical paradox, we move on to Mystery Two.
   If we have a brake orbit whose ends-- the two brake triangles - are related by   a  symmetry $F$ as above  we can use that symmetry to
  extract non-trivial information about the configuration $q(T/4)$ at the quarter period, $T/4$  being 
    half-way between the two brake times of $t =0$ and $t = T/2$.  
  Suppose then, that $F(q(0)) = q(T/2)$.  Consider the new solution
  $F(q(t))$ at time $t =0$ its initial condition is shared with that of $q(t)$ at
  time $t = T/2$: namely, it is brake with configuration $q(T/2)$.  But $q(t + T/2)$
  solve Newton's equations and has precisely these initial conditions.
  ( Like any  autonomous ODE (equation (\ref{N})) enjoys the time translational  symmetry: 
  If $q(t)$ solves so does $q(t + t_0)$ for any time $t_0$.)  It follows that
  \beq
  q(t + T/2) = F(q(t)).
  \label{eq: symm}
  \eeq  Now take $t = -T/4$ and use the time reversal symmetry to
  conclude that 
  \beq
 q(T/4) = F(q(T/4))
 \label{eq: quarter_symm}
 \eeq  The midway point must  a fixed point of our symmetry $F$!
  
  Return to figure \ref{fig18Gofen}. A reflection $F$ relates its two brake triangles. Reflections are symmetries of the three-body problem.   
  Take $F$ to be this reflection and $\ell$ its line of  reflection.   $F$'s fixed points are 
  collinear ``triangles'' in which all three masses lie on $\ell$.   In this way we have solved the mystery around that figure:  why all three masses form a syzygy at the mid time, and indeed at the same time a moment's contemplation of equation (\ref{eq: symm}) shows that we have also 
   accounted for the overall reflectional symmetry
  of that orbit as being a consequence of the symmetric relation between its endpoints.  
   
 In figure \ref{fig14Gofen} the two brake triangles are related by $F = R \circ \sigma$ where $R$ is
 rotation by 180 degrees and $\sigma$ interchanges two of the masses.  The fixed point set of such an $F$ is the
 set of ``Euler configurations'' :  degenerate collinear triangles with the non-interchanged mass forming the midpoint of the other two.
 This fact, and of course equation (\ref{eq: symm})  matches Gofen's data.

 
 \section{A hole in shape space and  harmonic oscillators}  
 Our mathematical hero Arnol'd had a saying he was fond of sprinkling into his lectures which  was a variation  of
 the phrase  ``the exception that proves the rule''.
 We present our exception.  
 
 Observe that we can  run our shape space argument   for any potential invariant under isometries.  
 One such potential  is the harmonic oscillator potential $V = \Sigma k_{ab} r_{ab}^2$ with $k_{ab} > 0$ spring constants.
 The act of replacing  Newton's potential (equation (\ref{eq: potential})) by this quadratic potential    corresponds to replacing the gravitational force by Hooke's spring forces.
 The corresponding ODEs are linear of the   form   $\ddot q  = - A q $ where $A$
 is a matrix depending on the masses and springs and which is positive definite on the center-of-mass subspace $\Sigma m_a q_a = 0$
 and which leaves  this subspace invariant.   Choose an eigenbasis $E_i$ for $A$ restricted to this subspace
 so that $A E_i = - \omega_i ^2 E_i$ with $\omega_i ^2 > 0$ the eigenvalues.  Then $q(t) =  \cos(\omega_i  t) E_i$ is a brake solution shuttling back and forth
 between $E_i$ and $-E_i$.  But $-E_i$ corresponds to rotating $E_i$ by $180$ degrees.  Pushed down to shape space,
 this brake solution connects the shape corresponding to $E_i$ to itself, contradicting my  alleged 
 ``shape-space thinking'' theorem that such a solution is  impossible. This is the exception that proves the rule. 
 
 Another paradox.  What's happening? Is our theorem true or not?    
 
  The resolution of this apparent paradox involves the projection of $0$ to shape space, 
 $0$ representing total collision.  The map from configuration space to shape space, upon  restriction  to the center-of-mass subspace,
 fails to be a submersion exactly at $0$.    We should view the shape $0$ of total collision as a singularity in shape space.
 (Indeed, for $N > 3$ it is a topological singularity since the shape space for the planar N-body problem is the cone over complex projective
 space of complex dimension $N-2$.)  The relation between the dynamics upstairs and downstairs breaks down at total collision. 
 Our eigenvector-based
 solution   above  passes through $0$   at time   $t = \pi/(2 \omega_i)$,
 and our shape-space argument fails for solutions   passing through $0$.     
 
 We can derive an  alternative   resolution to this paradox   by following the implications of equation
 (\ref{eq: symm}).  A rotation $F = R$ is a symmetry and so that equation  hold for any  brake solution
 whose  end triangles are related by a rotation $R$.  Then equation (\ref{eq: quarter_symm})
 asserts that $q(T/4)$ is a fixed point of the rotation.  But the only centered configuration
 invariant under a nontrivial rotation is the zero configuration $0 = (0,0,0)$, the configuration  representing total collision.  Our  periodic  brake  solution  
 must pass through   total collision half-way between its two ends!   That's fine for
 the harmonic oscillator.   No problem.  For the gravitational N-body problem  total collision
 acts like an essential singularity - a hole in shape space if you will -  through which  there is no consistent way to travel beyond and we have to stop
 the dynamics at total collision and call it quits.  So there's no such brake orbit for the planar N-body problem.
 
 More information can be extracted from equation (\ref{eq: symm}).   Evaluating the  equation  at $t = -T/2$ and $t =0$
 to get $q(0) = F(q(T/2))$ and   $q(T/2) = F(q(0))$ so that    
 $q(0) = F^2 (q(0))$.     If the triangle $q(0)$ is in general position,  or even if it is a degenerate collinear triangle
 but  $F$ is of the form $R$ or $R \circ \sigma$,
 then this fixed point relation  implies that  $F^2 = Id.$ (In other terms,
  the rotation group 
 acts freely on configuration space away from triple collision, so that $R^2 (q(0)) = q(0)$
 implies that $R^2 = Id.$) 
 When $F = R \circ \sigma$ we have  $F^2 = R^2$, so either way, when $F = R$ or $R \circ \sigma$,  we get   $R^2 = Id$.   
   {\it The only $R$'s that solve the identity $R^2 = Id$ are rotations by    180 degrees! }
 This $R$ is also known as  central  inversion:   $R q_a = -q_a$.     
 The original solution  Gofen showed me, figure \ref{fig14Gofen} above,
 has its brake triangles related by   central inversion, but again, related as  {\it unlabelled} triangles.  
Again, this matches Gofen's data:  all his nontrivial brake triangles, when related by symmetries,
are related by central inversion.

 I sure hope Gofen does not wander back   out into  the land of equal mass periodic 
 brake orbits this summer  and comes back to me  with a solution  whose brake triangles are related by an 
  $F$   of the form $R \sigma$  with the $R$    being a 45 degree rotation!   I will not know how to 
  resolve the resulting paradox. The remainder of my summer vacations with family would be threatened with ruin!
  Wait till next spring, please, Alex, for alerting me to  such a paradox. 
  
  Mystery one, the mystery of 12 out of 30 of these ``first''  equal mass collision-free periodic brake orbits  having extra symmetries,  
  remains a mystery.   
  
 \section{End Note. Gofen's  Taylor Center}

Alex Gofen, whose figures grace this paper,  asked me to say a few things
about the   ``Taylor Center'' that he runs and which generated three of the displayed figures. 
 
 The name Taylor Center  stands for two things:
\begin{itemize}
\item{}
a comprehensive resource dedicated to particular mathematical problems; and...
\item{}
 a software for Windows - the advanced ODE solver [b] based on the modern Taylor integration. This ODE solver offers several unique features: numerical and graphical.
\end{itemize}
As a numerical tool, it employs the most accurate Intel float point type called extended with 63-bit mantissa, integrating with order 30 or higher, and providing several methods of accuracy control up to all available 63 binary digits. 

As a graphical tool, it offers high resolution graphics, plotting the trajectories as real time animation: both in 2D and 3D stereo (viewable via red/blue glasses). Thanks to such graphics, this software may serve as Lab-works in various fields of applied mathematics.  A few of such lab topics were already posted [c]: for example, the three types of the rigid body motion; selected samples in celestial mechanics. The Lab-works library keeps growing.

[a] \url{http://taylorcenter.org/}

[b] \url{http://taylorcenter.org/Gofen/TaylorMethod.htm}

[c] \url{http://taylorcenter.org/Exploratorium/}

\date{today}

\end{document}